\documentclass[12pt,a4paper,reqno]{amsart}

\usepackage{amsmath}
\usepackage{amstext}
\usepackage{amsfonts}
\usepackage[mathscr]{eucal}
\usepackage{amscd}
\usepackage{latexsym}
\usepackage{amssymb}
\usepackage{diagrams}

\newtheorem{theorem}[]{Theorem}
\newtheorem{lemma}[theorem]{Lemma}
\newtheorem{proposition}[theorem]{Proposition}
\newtheorem{definition}[theorem]{Definition}

\theoremstyle{definition}

\newcommand{\zbb}{\mathbb{Z}}
\newcommand{\gbb}{\mathbb{G}}

\newcommand{\mi}{\text{-}}
\newcommand{\enpr}{\hfill $\Box $}
\newcommand{\supar}[1]{\overset{#1}\rightarrow}

\newcommand{\prf}{\noindent\textbf{Proof. }}
\newcommand{\rw}{\rightarrow}

\newcommand{\vep}{\varepsilon}

\newcommand{\id}{\mathrm{id}}
\newcommand{\der}{\mathrm{Der}}
\newcommand{\hm}{\mathrm{Hom}}

\newcommand{\ifr}{\mathfrak{I}}
\newcommand{\im}{\mathrm{Im}\,}

\newcommand{\hund}[1]{{\hm}_{#1}}
\newcommand{\otund}[1]{{\otimes}_{#1}}

\newcommand{\ab}{\mathbf{Ab}}
\newcommand{\cm}{\mathbf{CM}}
\newcommand{\tgmu}{(T,G,\mu)}
\newcommand{\abde}{(A,B,\delta)}

\newcommand{\ccl}{\mathcal{C}}
\newcommand{\gcl}{\mathcal{G}}
\newcommand{\remb}{\rightarrowtail}
\newcommand{\ronto}{\twoheadrightarrow}

\newcommand{\rt}{\rtimes}

\newcommand{\set}{\mathbf{Set}}

\newcommand{\ttil}{\widetilde{\times}}
\newarrow{Eq}{}{=}{=}{=}{}
\newarrow{LR}{<}{-}{-}{-}{>}
\newarrow{Dts}{}{.}{.}{.}{}
\newarrow{Embedline}{>}{-}{-}{-}{-}
\newarrow{Intoar}{C}{-}{-}{-}{>}

\small%
\normalsize

\begin{document}%

\makeatletter
\renewcommand{\descriptionlabel}[1]{\hspace\labelsep \upshape\bfseries #1}
\makeatother

\title[On the non--balanced property of the category of ..]{On the non--balanced property of the category of crossed modules in groups.}
\author{ Simona Paoli }
\email{simona@maths.warwick.ac.uk}
\address{Mathematics Institute\\ University of Warwick\\ Coventry CV4 7AL \\ UK
 }

\subjclass{18G50 (18C15)}

\keywords{Crossed module, category of interest.}

\begin{abstract}
An algebraic category $\mathcal{C}$ is called balanced if the cotriple cohomology of any
object of $\mathcal{C}$ vanishes in positive dimensions on injective coefficient modules.
Important examples of  balanced and of non-balanced categories occur in the literature.
In this paper we prove that the category of crossed modules in groups is non-balanced.
\end{abstract}

\maketitle

\section*{Introduction}\label{intro}
The (co)homology theory of algebraic objects was studied using cotriple resolutions since
the time of Beck's work \cite{beck}. In the category of groups cotriple (co)homology
recovers, up to a dimension shift, ordinary group (co)homology, and similar results hold
for Lie algebras and associative algebras over a field \cite{barr}.

An axiomatization of the Barr-Beck setting was given by Orzech \cite{orz} who defined
'categories of interest', which include the examples mentioned above.

The interpretation of the cohomology groups $H^{n}(G,A)$ of a group $G$ with coefficients
in a $G$-module $A$ was given independently by various authors (see MacLane's historical
note \cite{mac2}) in terms of equivalence classes of crossed $n$-fold extensions of $G$
by $A$.

In \cite{vale} Vale proved that cotriple cohomology in a category of interest can be
interpreted in terms of equivalence classes of crossed $n$-fold extensions provided that
the the cotriple cohomology vanishes in positive dimensions on injective coefficients
modules. A category of interest with this property is called, in this context, a balanced
category.

Examples of balanced categories include groups and Lie algebras.  The main example of a
non-balanced category is commutative algebras. The corresponding cotriple theory is in
this case Andr\'{e}--Quillen cohomology \cite{qui}.

Crossed modules in groups were introduced by Whitehead \cite{whi} as algebraic models of
connected 2-types. One of the first studies of purely algebraic aspects of crossed
modules in groups was the work of Norrie \cite{nor} which investigated actions internally
to the category  of crossed modules. In \cite{ccg} the authors studied crossed modules as
a category tripleable over the category of sets, and they introduced a corresponding
(co)homology theory with trivial coefficients. A more general class of local coefficients
was considered in \cite{pao1}.

The purpose of this paper is to answer the question whether the algebraic category of
crossed modules in groups is balanced. In Theorem \ref{cronon.the1} we prove that crossed
modules is non-balanced.

This paper is organized as follows. In Section 1 we provide some background. In
\ref{catint} we recall a criterion, due to Cegarra \cite{ceg}, which is satisfied in any
balanced category of interest. In \ref{interest} we recall how the category of crossed
modules can be considered  as a category of interest via its equivalence with the
category of cat$^{1}$-groups.

In the first part of Section 2  we characterize  modules in the category of crossed
modules in groups, and reconcile the notion of action  in the sense of categories of
interest with the one described by Norrie. In the second part of Section 2 we investigate
Cegarra's criterion in the category of crossed modules, and eventually prove that it does
not hold, showing that crossed modules is non-balanced.

\medskip
\section*{Acknowledgments}
 I would like to thank the Macquarie University Mathematics
Department for their hospitality during the writing of this paper, and the members of the
Australian Category Seminar for providing a stimulating working environment during my
visit. I am indebited in particular to Dr Steve Lack for showing me how to view
cat$^{1}$-groups as a variety of universal algebras. This paper is based on some results
from my PhD thesis \cite{pao2}; I would like to thank my PhD supervisor Dr Alan Robinson
for his support and encouragement.

\bigskip

\section{Preliminaries}\label{prel}
\subsection{Categories of interest and Cegarra's criterion.}\label{catint} Let $\ccl$ be
a category of interest in the sense of Orzech \cite{orz}. Recall that this consists of a
category satisfying the following axioms:
\begin{description}
  \item[1)]\hspace{7pt} There is a triple $(\top,\eta,\mu)$ on $\set$ such that $\top(\emptyset)=\{\cdot\}$ (a
  one point set) and $\ccl$ is equivalent to $\set^\top$.
\end{description}
Let $\set_*$ denote the category of pointed sets, with basepoint preserving maps. From
above, $\ccl$ is pointed with $(\top(\emptyset),\mu_\emptyset)$ as zero object. Also,
$\ccl$ is tripleable over $\set_*$.
\begin{description}
  \item[2)] \hspace{7pt} The underlying set functor $U:\ccl\rw \set_*$ factors through the
  category of groups.
  \item[3)] \hspace{7pt} All operations in $\ccl$ are finitary.
  \item[4)] \hspace{7pt} There is a generating set $\Omega$ for the operations in $\ccl$ and $\Omega=
  \Omega_0\cup \Omega_1\cup\Omega_2$, where $\Omega_i=$ set of $i$-ary operations in
  $\Omega$. Moreover $\Omega$ includes the identity, inverse and multiplication associated
  with the group structure; let $\Omega'_2=\Omega_2\verb"\"\{+\}$, $\;\Omega'_1=\Omega_1
  \verb"\"\{-\}$. We assume that if $*\in \Omega'_2$ then $*^o$ defined by $x*^oy=y*x$ is
  also in $\Omega'_2$.
  \item[5)] \hspace{7pt} For each $*\in\Omega'_2$, $\;a*(b+c)=a*b+a*c$.
  \item[6)] \hspace{7pt} For $\omega\in \Omega'_1$, $\;*\in\Omega'_2$, $\;\omega(a+b)=\omega(a)+
  \omega(b)$ and $\omega(a*b)=\omega(a)*b$.
  \item[7)] \hspace{7pt} For each $*\in\Omega'_2$, $\;a+(b*c)=(b*c)+a$.
  \item[8)] \hspace{7pt} For each ordered pair $(\cdot,*)\in\Omega'_2\times\Omega'_2$
  there is a word $w$ such that
\begin{equation*}
\begin{split}
  (x_1\cdot x_2)*x_3= & w(x_1(x_2x_3),x_1(x_3x_2),(x_2x_3)x_1,(x_3x_2)x_1,\\
  & x_2(x_1x_3),x_2(x_3x_1),(x_1x_3)x_2,(x_3x_1)x_2),
  \end{split}
\end{equation*}
where juxtaposition represents an operation in $\Omega'_2$.
\end{description}
  We recall the notion of module in a
category of interest. An object $A$ of $\ccl$ is called singular if $A$ is an abelian
group and if $a_1*a_2=0$ for all $a_1,a_2\in A$, $\:*\in \Omega_2\verb"\"\{+\}$. Given an
object $R$ of $\ccl$, an $R$-\emph{module} consists of a singular object $A$ and of a
split short exact sequence in $\ccl$
\begin{equation*}
  A\remb E\;\pile{\ronto \\ \leftarrow}\; R.
\end{equation*}
It follows that $E\cong A\ttil R$, where $A\ttil R=A\times R$ as a set, with operations
\begin{align*}
   & \omega(a,r)=(\omega(a),\omega(r)),\quad\text{for each}\;\omega\in \Omega_1 \symbol{92} \{-\}, \\
   & (a',r')+(a,r)=(a'+s(r')+a-s(r'),r'+r),\\
   & (a',r')*(a,r)=(a'*a+a'*s(r)+s(r')*a,r'*r), \quad\text{for each}\; * \in \Omega_2 \symbol{92}\{+\}.
\end{align*}
A \emph{morphism of $R$-modules} consists of a morphism $f:A\rw A'$ in $\ccl$ inducing a
commutative diagram of split extensions
\begin{diagram}[s=2em]
    A & \rEmbed & A\ttil R & \rOnto & R\\
    \dTo^f && \dTo_{(f,\id)} && \dEq\\
    A' & \rEmbed & A'\ttil R' & \rOnto & R
\end{diagram}
The category $R$-Mod of $R$-modules in $\ccl$ is equivalent to the category
$(\ccl/R)_{ab}$ of abelian group objects in $\ccl/R$. Further, it can be proved
\cite[p.30]{ceg} that $R$-Mod has enough injectives.

Given an object $Y\supar{f} R$ in the slice category $\ccl/R$, any $R$-module $A$ is also
a $Y$-module via the map $f$. There is a derivation functor $\der(\mi,A):\ccl/R\rw \ab$,
such that
\begin{equation}\label{catint.1}
  \der(Y,A)\cong\hund{\ccl/R}(Y,A\ttil R).
\end{equation}
An explicit description of a derivation $D:Y\rw A$ in $\ccl$ can be given, see for
instance \cite[p.36]{ceg}. The forgetful functor $(\ccl/R)_{ab}\rw\ccl/R$ has a left
adjoint $D_R:\ccl/R\rw(\ccl/R)_{ab}$. Since $(\ccl/R)_{ab}$ is equivalent to $R$-Mod, and
$A\ttil R\rw R$ corresponds to $A$ under this equivalence, by adjointness and by
(\ref{catint.1}) it follows that
\begin{equation}\label{catint.2}
  \der(Y,A)\cong\hund{R\mi\mathrm{Mod}}(D_R(Y),A).
\end{equation}
Let $\gbb$ be the cotriple in $\ccl$ arising from the forgetful functor $\ccl\rw\set$ and
its left adjoint. This cotriple induces a cotriple on the slice category $\mathcal{C}/R$
which we still denote by $\gbb$. The $n^{th}$ cotriple cohomology of $R$ with
coefficients in the $R$-module $A$ is defined by
\begin{equation*}
  D^n(R\;,A)=H^{n}\der(\gbb_{*}R,A)
\end{equation*}
where $\gbb_{*}R\rw R$ is the cotriple resolution. A category of interest $\ccl$ is said
to be \emph{balanced} if for every object $R$ of $\ccl$, $D^n(R,I)=0$ for each $n>0$,
whenever $I$ is an injective $R$-module.

Let $N\remb Y \overset{f}{\ronto}R$ be a short exact sequence in $\ccl$. The commutator
subobject $[N,N]$ is the ideal of $N$ generated by the elements
$\{n_1+n_2-n_1-n_2,n_1*n_2\;|\;n_1,n_2\in N,\;*\in \Omega_2\symbol{92}\{+\}\}$; the
quotient $N/[N,N]$ is an $R$-module, with the action of $R$ on $N/[N,N]$ given by
\begin{equation}\label{catint.3}
  r\cdot[n]=[y+n-y],\qquad r*[n]=[r*n]
\end{equation}
where $f(y)=r$, $y\in Y$, $r\in R$, $n\in N$.

\begin{theorem}{\cite[(1.7.11)]{ceg}}\label{catint.the1}
Let $\ccl$ be a category of interest, $N\remb Y\ronto R$ a short exact sequence in
$\ccl$. Then there exists a natural exact sequence of $R$-modules
\begin{equation*}
  \frac{N}{[N,N]}\supar{j}D_R(Y)\ronto D_R(R).
\end{equation*}
Moreover, if $\ccl$ is balanced, then $j$ is a monomorphism.
\end{theorem}
\bigskip

\noindent\textbf{Examples:}
\begin{description}
  \item[a)] Let $\ccl$ be the category of groups, and let $N\remb G\ronto G'$ be a short
  exact sequence in groups. The corresponding 3-term exact sequence is
  $N_{ab}\remb\zbb G'\otund{\zbb G}\ifr_{G}\ronto\ifr_{G'}$. Given a $G$-module
  $A$, $D^n(G,A)\cong H^{n+1}(G,A)$ for $n>0$ \cite{beck}, hence $D^n(G,I)=0$ for
  $n>0$ and $I$ an injective $G$-module. The category of groups is a balanced category.
  \item[b)] Let $\ccl$ be the category of commutative algebras. Given the short exact
  sequence of commutative algebras $N\remb R\ronto R'$, the corresponding 3-term exact sequence is
  $N/N^2\supar{j}R\otund{A}\Omega_{R'}\ronto \Omega_{R'}$, where $\Omega_R$ is the module
  of K\"{a}hler differentials. The map $j$ is in general not injective, hence commutative
  algebras is an example of a non-balanced category.
\end{description}
\subsection{Crossed modules as a category of interest.}\label{interest}
Recall that the category $\cm$ of crossed modules in groups has objects the triples
$\tgmu$ where $\mu:T\rw G$ is a group homomorphism, $G$ acts on $T$ and for all $t,t'\in
T,\; g\in G$,
\begin{equation*}
  \mu(\,^{g}t)=g\mu(t)g^{-1}, \qquad ^{\mu(t)}t'=tt't^{-1}.
\end{equation*}
A morphism of crossed modules is a pair of group homomorphisms $(f,h):\tgmu\rw
(T',G',\mu')$ such that $\mu'f=h\mu$ and $f(\,^{g}t)=\:^{h(g)}f(t)$, $t\in T$, $g\in G$.
The category $\cm$ is equivalent to the category $\ccl^1\gcl$ of cat$^1$-groups
\cite{lod}. Objects of $\ccl^1\gcl$ are triples $(G,d_0,d_1)$ where $d_0,d_1:G\rw G$ are
group homomorphisms and
\begin{equation}\label{interest.1}
  d_0d_1=d_1,\qquad d_1d_0=d_0,\qquad [\ker d_0,\ker d_1]=1.
\end{equation}
A morphism of cat$^1$-groups $f:(G,d_0,d_1)\rw (G',d'_0,d'_1)$ is a group homomorphism
$f:G\rw G'$ such that $fd_i=d'_if,\;i=0,1$. Since $\ker d_i=\{d_i(x)x^{-1}\;|\;x\in G\}$,
$i=0,1$ the identities (\ref{interest.1}) are equivalent to
\begin{equation*}
  d_0d_1=d_1,\qquad d_1d_0=d_0,\qquad d_0(x)x^{-1}d_1(y)y^{-1}=d_1(y)y^{-1}d_0(x)x^{-1},
  \; x,y\in G.
\end{equation*}
It follows that $\ccl^1\gcl$ is a category of universal algebras. The generating set of
operations is $\Omega=\Omega_0\cup\Omega_1\cup\Omega_2$, $\Omega_0=\{0\}$,
$\Omega_1=\{-\}\cup \{d_0,d_1\}$, $\Omega_2=\{+\}$, where $0,-,+$ denote group identity,
inverse and multiplication and the following identities hold; for all $x,y\in G$,
\begin{align*}
   & d_0(x+y)=d_0(x)+d_0(y),\qquad d_1(x+y)=d_1(x)+d_1(y), \\
   & d_0d_1(x)=d_1(x),\qquad d_1d_0(x)=d_0(x),\\
   & d_0(x)x^{-1}d_1(y)y^{-1}=d_1(y)y^{-1}d_0(x)x^{-1}.
\end{align*}
It follows that the forgetful functor $U:\ccl^1\gcl\rw\set$ has a left adjoint
\cite{mac}. An explicit and very useful description of the left adjoint to $U$ was given
in \cite{ccg}.

The above description of $\ccl^1\gcl$ also makes this category into a `category of
interest' in the sense of Orzech \cite{orz}. Since $\cm$ is equivalent to $\ccl^1\gcl$,
we can therefore consider $\cm$ itself as a category of interest.

\section{The non-balanced property of crossed modules.}
\subsection{The notion of module in the category $\cm$.}\label{modcm}
\begin{definition}\label{modcm.def1}
    Let $\tgmu$ and $\abde$ be crossed modules, and let $\gcl_{\tgmu}$, $\gcl_{\abde}$ be
    the corresponding cat$^1$-groups. We say that $\abde$ is a $\tgmu$-module if
    $\gcl_{\abde}$is a $\gcl_{\tgmu}$-module in the sense of categories of interest.
\end{definition}
\begin{lemma}\label{modcm.lem1}
    Given crossed modules $\tgmu$ and $\abde$, the following are equivalent:
\begin{description}
  \item[i)] $\abde$ is a $\tgmu$-module.
  \item[ii)] $\abde$ is an abelian crossed module and there is a split short exact
  sequence in $\cm$
\begin{equation}\label{modcm.1}
  \abde \overset{i}{\remb}(T',G',\mu')\;\pile{\ronto\\ \underset{s}{\leftarrow}}\;\tgmu
\end{equation}
  \item[iii)] $\abde$ is an abelian crossed module and there is an action of $\tgmu$ on
  $\abde$in the sense of \cite{nor}, that is there is a crossed module morphism $(\vep,\rho):\abde\rw
  \mathrm{Act}\tgmu$.
\end{description}
\end{lemma}
\prf

\bigskip

\noindent i) $\Leftrightarrow$ ii)

By definition $\abde$ is a $\tgmu$-module if and only if $\gcl_{\abde}$ is a singular
object in $\ccl^1\gcl$ and there is a split short exact sequence in $\ccl^1\gcl$ of the
form
\begin{equation}\label{modcm.2}
    \gcl_{\abde}\remb \gcl \;\pile{\ronto\\ \leftarrow}\;\gcl_{\tgmu}.
\end{equation}
Here $\gcl_{\abde}\!=\!(A\rt B,d_0,d_1)$, $d_0(a,b)=(1,b)$, $ d_1(a,b)=(1,\delta(a)b)$.
By definition, $\gcl_{\abde}$ is a singular object in $\ccl^1\gcl$ if and only if $A\rt
B$ is an abelian group. This is equivalent to $\abde$ being an abelian crossed module. In
fact, if $A\rt B$ is abelian, both $A$ and $B$ are abelian groups and for all $a\in A, \;
b\in B $ $(a,b)=(a,0)+(0,b)=(0,b)+(a,0)=(\:^{b}a,b)$ so that $a=\:^{b}a$. Hence $\abde$
is an abelian crossed module. Conversely, if $\abde$ is an abelian crossed module, then
$A\rt B=A\oplus B$ is an abelian group.

Let $(T',G',\mu')$ be the crossed module corresponding to the cat$^1$-group $\gcl$. Then
(\ref{modcm.2}) holds if and only if there exists a split short exact sequence in $\cm$
of the form
\begin{equation}\label{modcm.3}
  \abde\remb(T',G',\mu')\;\pile{\ronto\\ \leftarrow}\;\tgmu.
\end{equation}
In conclusion i) is equivalent to the existence of a short exact sequence (\ref{modcm.2})
with $\gcl_{\abde}$ a singular cat$^1$-group; in turn this is equivalent to the short
exact sequence (\ref{modcm.3}) with $\abde$ an abelian crossed module, which is  ii).

\bigskip

\noindent ii) $\Leftrightarrow$ iii)

If ii) holds, then the split short exact sequence (\ref{modcm.1}) induces split short
exact sequences of groups $A\remb T'\;\pile{\ronto\\ \leftarrow}\;T$ and $B\remb
G'\;\pile{\ronto\\ \leftarrow}\;G$, so that there are induced actions of $T$ on $A$ and
of $G$ on $B$, and $T'=A\rt T$, $G'\cong B\rt G$. Since the maps $i$ and $s$ in
(\ref{modcm.1}) are injective, we can identify $\tgmu$ with $s\tgmu$ and $\abde$ with
$i\abde$. Hence, by (\ref{modcm.1}), we can regard $\abde$ and $\tgmu$ as subcrossed
modules of $(T',G',\mu')$. With these identifications, we have
\begin{description}
  \item[a)] $\abde$ is a normal subcrossed module of $(T',G',\mu')$.
  \item[b)] $T'=AT,\quad G'=BG$.
  \item[c)] $A\cap T=1, \quad B\cap G=1$.
\end{description}
By \cite[p. 135]{nor} it follows that there is a morphism of crossed modules
$(\vep,\rho):\abde\rw\mathrm{Act}\tgmu$, $\vep_a(g)=a-\:^{g}a$,
$\rho(b)=(\alpha_b,\gamma_b)$, $\alpha_b(t)=\:^{b}t$, $\gamma_b(t)=btb^{-1}$. Hence ii)
implies iii).

Conversely, if \: iii) holds, by \cite{nor} there exists a split short exact sequence of
crossed modules
\begin{equation*}
  \abde\remb(A\rt T, B\rt G, (\delta,\mu))\;\pile{\ronto\\ \leftarrow}\;\tgmu
\end{equation*}
where the crossed module action of $B\rt G$ on $A\rt T$ is given by
\begin{equation*}
  ^{(b,g)}(a,t)=(\:^{b}(\:^{g}a)-\vep(\:^{g}t)(b),\:^{g}t)=(\:^{g}a-\vep(\:^{g}t)(b),\:^{g}t).
\end{equation*}
Hence ii) holds.\enpr
\begin{definition}\label{modcm.def2}
    Let $\abde$, $(A',B',\delta')$ be $\tgmu$-modules. A morphism
    $f:\abde\rw(A',B',\delta')$ of $\tgmu$-modules is a crossed module morphism such that
    the corresponding morphism of cat$^1$-groups
    $f:\gcl_{\abde}\rw\gcl_{(A',B',\delta')}$ is a morphism of $\gcl_{\abde}$-modules in
    the sense of categories of interest.
\end{definition}
\begin{lemma}\label{modcm.lem2}
    A morphism of $\tgmu$-modules $(r,s):\abde\rw(A',B',\delta')$ consists of a pair of
    morphisms of abelian groups such that
\begin{equation*}
  \delta' r=s\delta, \qquad
  r(\:^{g}a)-r(\vep(\:^{g}t)(b))=\:^{g}r(a)-\vep'(\:^{g}t)(s(b))
\end{equation*}
for $a\in A$, $b\in B$, $g\in G$, $t\in T$, where $(\vep,\rho):\abde\rw\mathrm{Act}\tgmu$
and $(\vep',\rho'):(A',B',\delta')\rw\mathrm{Act}\tgmu$.
\end{lemma}
\prf By Definition \ref{modcm.def2} and Lemma \ref{modcm.lem2}, $(r,s)$ is a morphism of
$\tgmu$-modules if and only if there is a commutative diagram of split short exact
sequences
\begin{diagram}[s=2em]
    \abde & \rEmbed & (A\rt T, B\rt G, (\delta,\mu)) & \pile{\rOnto\\ \lTo} & \tgmu\\
    \dTo^{(r,s)} && \dTo_{((r,\id_T),(s,\id_G))} && \dEq\\
    (A',B',\delta') & \rEmbed & (A'\rt T, B'\rt G, (\delta',\mu')) & \pile{\rOnto\\ \lTo} & \tgmu
\end{diagram}
in $\cm$. This is equivalent to $(r,s)$ and $((r,\id_T),(s,\id_G))$ being crossed module
morphisms. Since $\abde$ and $(A',B',\delta')$ are abelian crossed modules, $(r,s)$ is a
crossed module morphism if and only if $\delta'r=s\delta$. For $((r,\id_T),(s,\id_G))$ to
be a crossed module morphism we further require that
\begin{equation*}
  (r,\id_T)(\:^{(b,g)}(a,t))=\:^{(s(b),g)}(r(a),t)
\end{equation*}
for all $a\in A$, $b\in B$, $t\in T$, $g\in G$. Hence
\begin{equation*}
  (r,\id_T)(\:^{g}a-\vep(\:^{g}t)(b),\:^{g}t)=(\:^{g}r(a)-
  \vep'(\:^{g}t)(s(b)),\:^{g}t)\; ;
\end{equation*}
that is
\begin{equation*}
  r(\:^{g}a)-r(\vep(\:^{g}t)(b))=\:^{g}r(a)-\vep'(\:^{g}t)(s(b)).
\end{equation*}
\enpr
\subsection{Crossed modules is non-balanced.}\label{cronon}
\begin{lemma}\label{cronon.lem1}
    Let $M$ be an abelian group and consider an extension of crossed modules
\begin{equation*}
  (N,G,\nu)\remb\tgmu\overset{(f,0)}{\ronto}(M,1,0).
\end{equation*}
Let $J=[G,N][T,T]$. Then $(T/J, G_{ab},\overline{\mu})$ is a $(M,1,0)$-module with action
$(\vep',1):(M,1,0)\rw \mathrm{Act}(T/J,G_{ab},\overline{\mu})$,
$\overline{\mu}[t]=[\mu(t)]$, $\vep'(m) [g]=[t]-[\,^{g}t]$, $m=f(t)$, $t\in T$, $g\in G$.
Moreover,
\begin{equation}\label{cronon.1}
  \der(\tgmu,\abde)\cong\hund{(M,1,0)\mi\mathrm{Mod}}((T/J,G_{ab},\overline{\mu}),\abde).
\end{equation}
\end{lemma}
\prf We first check that
\begin{equation*}
(\vep',1):(M,1,0)\rw \mathrm{Act}(T/J,G_{ab},\overline{\mu})
\end{equation*}
 is a crossed module morphism. We recall from \cite{nor} that
\begin{equation*}
\mathrm{Act}(T/J,G_{ab},\overline{\mu})=(D(G_{ab},T/J)\;,
\mathrm{Aut}(T/J,G_{ab},\overline{\mu})(\theta,\sigma)),
\end{equation*}
where $D(G_{ab},T/J)$ is the group of Whitehead derivations.

The map $\vep'$ is well defined. In fact, if $m=f(t_1)=f(t_2)$, then $t_1t_2^{-1}\in N$
so that $[\,^{g}(t_1t_2^{-1})]=[t_1t_2^{-1}]$, $t_1t_2\in T,\;g\in G$. Hence
\begin{equation}\label{cronon.3}
  \vep'(f(t_1))[g]=[t_1]-[\,^{g}t_1]=[t_2]-[\,^{g}t_2]=\vep'(f(t_2)) [g].
\end{equation}
Choosing $t_1=t$, $t_2=\:^{g^{-1}}t$ in (\ref{cronon.3}), we obtain
$[\,^{g}t]-[t]+[\,^{g^{-1}}t]-[t]=0$ for all $t\in T,\;g\in G$. Hence, for all
$g_1,g_2\in G, \;t\in T$
\begin{equation*}
\begin{split}
   & [\,^{g_1g_2g_1^{-1}g_2^{-1}}t]-[t]=[\,^{g_1g_2g_1^{-1}g_2^{-1}}t]
   -[\,^{g_2g_1^{-1}g_2^{-1}}t]+\\
   & +[\,^{g_2g_1^{-1}g_2^{-1}}t]-[\,^{g_1^{-1}g_2^{-1}}t]+[\,^{g_1^{-1}g_2^{-1}}t]-[\,^{g_2^{-1}}t]+
   [\,^{g_2^{-1}}t]-[t]=\\
   & =[\,^{g_1}t]-[t]+[\,^{g_2}t]-[t]+[\,^{g_1^{-1}}t]-[t]+[\,^{g_2^{-1}}t]-[t]=0.
\end{split}
\end{equation*}
This proves that $[\,^{x}t]=[t]$ for all $x\in [G,G],\;t\in T$. Hence, if $[g_1]=[g_2]$,
$\vep'(m) [g_1]=\vep'(m) [g_2]$.

Notice that $\vep'(m)\in D(G_{ab},T/J)$. In fact, $\im
\vep'(m)\subseteq\ker\overline{\mu}$, since $\overline{\mu}(
[t]-[\,^{g}t])=[\mu(t)]-[g\mu(t)g^{-1}]=0$. It follows easily from \cite{nor} that
$\vep'(m)$ is a unit in $\der(G_{ab},T/J)$, hence is in $D(G_{ab},T/J)$. By \cite{nor}
for each $t,t'\in T, \;g\in G, \; m=f(t)$,
\begin{equation*}
\begin{split}
   & \theta\vep'(m)[t']=\vep'(f(t))\overline{\mu}[t']+[t']=[t]-[\,^{\mu(t')}t]+[t']=[t]-[t'tt'^{-1}]+
   [t']=[t'], \\
   & \sigma\vep'(f(t))[g]=\overline{\mu}\vep'(f(t))[g]+[g]=\overline{\mu}([t]-[\,^{g}t])+[g]=\![g].
\end{split}
\end{equation*}
Hence $(\theta,\sigma)\epsilon'=\mathrm{id}$, so that $(\epsilon',1)$ is a crossed module
morphism.

 By Lemma \ref{modcm.lem1}, a $(M,1,0)$-module consists of an abelian crossed
module $\abde$ and of a crossed module map $(\vep,1):(M,1,0)\rw
\mathrm{Act}\abde=(D(B,A),\mathrm{Aut}\abde,(\vartheta,\sigma))$. The corresponding split
extension in $\cm$ of $(M,1,0)$ by $\abde$ has the form
\begin{equation*}
  \abde\remb(A\oplus M, B,\overline{\delta})\;\pile{\ronto\\ \leftarrow}\;(M,1,0)
\end{equation*}
$\overline{\delta}(a,m)=\delta(a)$, $^{b}(a,m)=(a-\vep(m)(b),m)$, $a\in A$, $b\in B $,
$m\in M$.

By Lemma \ref{modcm.lem2}, $(r,s):\abde\rw(A',B',\delta')$ is a morphism
of $(M,1,0)$-modules if and only if $r,s$ are homomorphisms of abelian groups and
\begin{equation}\label{cronon.2}
  \delta'r=s\delta, \qquad r(\vep(m)(b))=\vep'(m)(s(b)),\qquad b\in B, \;m\in M,
\end{equation}
where $\vep$ and $\vep'$ are as in Lemma \ref{modcm.lem2}. By (\ref{catint.1}),
\begin{equation}\label{cronon.4a}
\der(\tgmu,\abde)\cong\hund{\cm/(M,1,0)}(\tgmu,(A\oplus M,B,\overline{\delta})).
\end{equation}
 We aim to show that
\begin{equation}\label{cronon.4}
 \begin{split}
 & \hund{\cm/(M,1,0)}(\tgmu,(A\oplus M,B,\overline{\delta}))\cong \\
 & \cong\hund{(M,1,0)\mi\mathrm{Mod}}((T/J,G_{ab},\overline{\mu}),\abde).
 \end{split}
\end{equation}
Let take $((D_1,f),D_2)\in \hund{\cm/(M,1,0)}(\tgmu,(A\oplus M,B,\overline{\delta}))$.
Then $D_1\in \der(T,A)$, $D_2\in \der(G,B)$ and for all $t\in T,\;g\in G$
\begin{equation}\label{cronon.5}
  D_2\mu=\delta D_1,\qquad D_1(\,^{g}t)=D_1(t)-\vep(f(t))(D_2(g)).
\end{equation}
Let $\varphi((D_1,f),D_2)=(\nu_1,\nu_2)$ where $\nu_1\in\hund{\zbb}(T/J,A)$,
$\nu_2\in\hund{\zbb}(G_{ab},B)$ are defined by $\nu_1[t]=D_1(t)$, $\nu_2[g]=D_2(g)$,
$t\in T,\;g\in G$.
\newline Since $\tgmu$ acts on $\abde$ via $(f,0)$ and $M$ acts trivially on $A$, then $T$ acts
trivially on $A$ and $G$ acts trivially on $B$. Hence $D_1(x)=0$ for all $x\in [T,T]$ and
$D_2(y)=0$ for all $y\in [G,G]$. Moreover from (\ref{cronon.5}) $D_1(\,^{g}n)=D_1(n)$ for
all $n\in N,\;g\in G$, so that $D_1(x)=0$ for all $x\in J$. It follows easily that
$\nu_1,\nu_2$ are well defined. By (\ref{cronon.5}),
\begin{equation*}
\begin{split}
   & \nu_2\overline{\mu}[t]=\nu_2[\mu(t)]=D_2\mu(t)=\delta D_1(t)=\delta\nu_1[t] \\
   & \nu_1[\,^{g}t]=\nu_1[t]-\vep(f(t))(\nu_2[g]).
\end{split}
\end{equation*}
From the definition of $\vep'$, we conclude that
\begin{equation}\label{cronon.6}
  \nu_2\overline{\mu}=\delta\nu_1,\qquad \nu_1(\vep'(m))[g]=\vep(m)(\nu_2[g]), \qquad
  m=f(t).
\end{equation}
Hence (\ref{cronon.2}) holds and $(\nu_{1},\nu_{2})$ is a morphism of $(M,1,0)$-modules.

Conversely, if $(\nu_1,\nu_2)\in
\hund{(M,1,0)\mi\mathrm{Mod}}((T/J,G_{ab},\overline{\mu}),\abde)$ then (\ref{cronon.6})
holds. Let $\psi(\nu_1,\nu_2)=((D_1,f),D_2)$ where $D_1(t)=\nu_1[t]$, $D_2(g)=\nu_2[g]$.
Since $T$ acts trivially on $A$ and $B$ acts trivially on $G$, $D_1\in\der(T,A)$ and
$B\in \der(G,B)$. From (\ref{cronon.6}), $D_2\mu=\delta D_1$ while taking $m=f(t)$ and
using the definition of $\vep'$ we obtain from the second identity in (\ref{cronon.6})
\begin{equation*}
  D_1(\,^{g}t)=D_{1}(\,t)-\vep(f(t))(D_2(g)).
\end{equation*}
Hence $((D_1,f),D_2)\in \hund{\cm/(M,1,0)}(\tgmu,(A\oplus M,B,\overline{\delta}))$.

It is straightforward to check that $\psi\varphi=\id$ and $\varphi\psi=\id$, proving the
isomorphism (\ref{cronon.4}). By (\ref{cronon.4a}), (\ref{cronon.1}) follows. \enpr

\bigskip

By the previous lemma, the left adjoint to the forgetful functor
$(\cm/(M,1,0))_{ab}\newline\rw\cm/(M,1,0)$ is given by
\begin{equation*}
  D_{(M,1,0)}\tgmu=(T/J,G_{ab},\overline{\mu})
\end{equation*}
where $\tgmu\rw(M,1,0)$ is an object of $\cm/(M,1,0)$ and $J$ is as in Lemma
(\ref{cronon.lem1}).
\begin{proposition}\label{cronon.pro1}
    Given the short exact sequence of crossed modules
\begin{equation*}
  (N,G,\nu)\remb\tgmu\overset{(f,0)}{\ronto}(M,1,0)
\end{equation*}
there is an exact sequence of $(M,1,0)$-modules
\begin{equation*}
  (N/[G,N],G_{ab},\overline{\nu})\overset{(u,\id)}{\rw}(T/J,G_{ab},\overline{\mu})
  \overset{(\overline{f},0)}{\ronto}(M,1,0)
\end{equation*}
where $J=[G,N][T,T]$, $u(n[G,N])=[n]$, $\overline{f}([t])=f(t)$, $n\in N$, $t\in T$. The
$(M,1,0)$-module structure of $(T/J,G_{ab},\overline{\mu})$ is as in Lemma
\ref{cronon.lem1}, and the $(M,1,0)$-module structure of
$(N/[G,N],G_{ab},\overline{\nu})$ is given by
$(\vep'',1):(M,1,0)\rw\mathrm{Act}(T/J,G_{ab},\overline{\mu})$, where
$\vep''(m)[g]=[t\:^{g}t^{-1}]$, $m=f(t)$, $g\in G$, $t\in T$.
\end{proposition}
\prf Let $\gcl_{(N,G,\nu)}$, $\gcl_{\tgmu}$, $\gcl_{(M,1,0)}$ be the cat$^1$-groups
corresponding to the crossed modules $(N,G,\nu)$, $\tgmu$, $(M,1,0)$ respectively. By
Theorem \ref{catint.the1}, there is an exact sequence of $\gcl_{(M,1,0)}$-modules
\begin{equation*}
  \gcl_{(N,G,\nu)}/[\gcl_{(N,G,\nu)},\gcl_{(N,G,\nu)}]\rw D_{\gcl_{(M,1,0)}}\gcl_{\tgmu}
  \ronto D_{\gcl_{(M,1,0)}}\gcl_{(M,1,0)}.
\end{equation*}
In the equivalent category of $(M,1,0)$-modules, $D_{\gcl_{(M,1,0)}}\gcl_{\tgmu}$
corresponds to $D_{(M,1,0)}\;{\tgmu}\,=\, (T/J,G_{ab},\overline{\mu})\:$ and
$D_{\gcl_{(M,1,0)}}\:\gcl_{(M,1,0)}$ corresponds to $D_{(M,1,0)}(M,1,0)= (M,1,0)$.

Since $\gcl_{(N,G,\nu)}=(N\rt G,s_0,s_1)$, $s_0(n,g)=(1,g)$, $s_1(n,g)=(1,\nu(n)g)$, it
is
\begin{equation*}
  \gcl_{(N,G,\nu)}/[\gcl_{(N,G,\nu)},\gcl_{(N,G,\nu)}]\cong ((N\rt G)_{ab},
  \overline{s}_0,\overline{s}_1).
\end{equation*}
where $\overline{s}_{i}$ is induced by $s_{i}\;, i=0,1.$  On the other hand it is not
hard to check that the map
\begin{equation}\label{cronon.7}
  \varphi:\bigg( \frac{N}{[G,N]}\rt\frac{G}{[G,G]}, u_0,u_1\bigg) \rw ((N\rt G)_{ab},
  \overline{s}_0,\overline{s}_1)
\end{equation}
$\varphi( [n],[g])=[(n,g)]$, $u_0( [n],[g])=(0,[g])$, $u_1( [n],[g])=(0,[\nu(n)g])$,
 is well defined and it is an isomorphism in
$\mathcal{C}^{1}\mathcal{G}$. Hence the crossed module corresponding to
$\gcl_{(N,G,\nu)}/[\gcl_{(N,G,\nu)},\gcl_{(N,G,\nu)}]$ is
$(N/[G,N],G_{ab},\overline{\nu})$.

By (\ref{catint.3}) the $\gcl_{(M,1,0)}$-module structure of
$\gcl_{(N,G,\nu)}/[\gcl_{(N,G,\nu)},\gcl_{(N,G,\nu)}]$ is given by
\begin{equation*}
  m\cdot[(n,g)]=[(t,1)(n,g)(t,1)^{-1}]=[(tn\:^{g}t^{-1},g)]
\end{equation*}
where $f(t)=m$, $[(n,g)]\in (N\rt G)_{ab}$, $t\in T$.

From the isomorphism (\ref{cronon.7}), the action of $\gcl_{(M,1,0)}$ on $\big(
\frac{N}{[G,N]}\rt {G}_{ab}, u_0,u_1\big)$ is therefore
\begin{equation}\label{cronon.8}
  m\cdot( [n],[g])=( [tn\:^{g}t^{-1}],[g])
\end{equation}
$m=f(t)$, $t\in T$, $g\in G$, $n\in N$. Consider the semidirect product in the category
of interest $\ccl^1\gcl$:
\begin{equation*}
  \frac{\gcl_{(N,G,\nu)}}{[\gcl_{(N,G,\nu)},\gcl_{(N,G,\nu)}]}\ttil\:
  \gcl_{(M,1,0)}\cong\bigg(\Big( \frac{N}{[G,N]}\rt G_{ab}\Big)\rt
  M,(u_0,0),(u_1,0)\bigg).
\end{equation*}
The crossed module corresponding to this semidirect product is
\begin{equation*}
(\ker(u_0,0),\im(u_1,0),(u_1,0)_{|\ker(u_0,0)})
\end{equation*}
with the action by conjugation of $\im(u_1,0)$ on $\ker(u_0,0)$. Hence we obtain from
(\ref{cronon.8}) the crossed module action
\begin{align*}
   & ^{((1,[g]),0)}(( [n],1),m)=((1,[g]),0)(( [n],1),m)((1,[g^{-1}]),0)= \\
   & =(( [n],[g]),m)((1,[g^{-1}]),0)=(( [n],[g])( [t\:^{g^{-1}}t^{-1}],[g^{-1}]),m)=\\
   & = (( [n]+[\,^{g}tt^{-1}],1),m).
\end{align*}
 There is clearly an isomorphism
\begin{equation*}
  (\ker(u_0,0),\im(u_1,0),(u_1,0)_{|\ker(u_0,0)})\cong \bigg( \frac{N}{[G,N]}\rt M,G_{ab},(
  \overline{\nu},0)\bigg).
\end{equation*}
Hence, by the previous calculation, the crossed module action of $G_{ab}$ on
$\frac{N}{[G,N]}\rt M$ is given by
\begin{equation*}
  ^{[g]}( [n],m)=( [n]+[\,^{g}tt^{-1}],m)
\end{equation*}
$g\in G$, $n\in N$, $m=f(t)$, $t\in T$. On the other hand, by \cite{nor}, the crossed
module action in the semidirect product crossed module $\big(
\frac{N}{[G,N]},G_{ab},\nu\big)\rt (M,1,0)\cong \big( \frac{N}{[G,N]}\rt M
,G_{ab},(\overline{\nu},0)\big)$ is given by
\begin{equation*}
  ^{[g]}( [n],m)=( [n]-\vep''(m) [g],m)
\end{equation*}
where $(\vep'',1):(M,1,0)\rw\mathrm{Act}(T/J,G_{ab},\overline{\mu})$. We deduce
$\vep''(m) [g]=-[\,^{g}tt^{-1}]=[t\,^{g}t^{-1}]$, $m=f(t)$, $g\in G$, $t\in T$. \enpr

\bigskip

\begin{theorem}\label{cronon.the1}
    The category of crossed modules in groups is non-balanced.
\end{theorem}
\prf We are going to show that there exists an extension of crossed modules such that the
corresponding 3-terms exact sequence is not short exact. By Theorem \ref{catint.the1}
this proves that crossed modules is non-balanced.

Let $M$ be an abelian group with $H_2(M)\neq 0$, and let $f:T\rw M$ be a surjection with
$T$ a free group. Let $N=\ker f$. Consider the extension of crossed modules
\begin{equation}\label{cronon.9}
  (N,T,i)\remb(T,T,\id)\overset{(f,0)}{\ronto}(M,1,0)
\end{equation}
where $i$ denotes the inclusion. The well known exact sequence $H_2(T)\rw H_2(M)\rw
N/[T,N]\rw T_{ab} \ronto M_{ab}$ reduces to
\begin{equation*}
  H_2(M)\remb N/[T,N]\supar{j}T_{ab}\ronto M
\end{equation*}
where $j(n[T,N])=n[T,T]$, $n\in N$. Since $H_2(M)\neq 0$, the map $j$ is not injective.
By Proposition \ref{cronon.pro1}, the 3-terms exact sequence associated to the extension
(\ref{cronon.9}) is
\begin{equation*}
  (N/[T,N],T_{ab},i)\supar{(j,\id)}(T_{ab},T_{ab},\id)\overset{(
  \overline{f},0)}{\ronto}(M,1,0)
\end{equation*}
and is not short exact. \enpr

\renewcommand{\refname}{References}

\bigskip

\end{document}